\documentclass[12pt]{article}

\usepackage{amsmath}
\usepackage{amssymb}
\usepackage[T1]{fontenc} 

\usepackage{colortbl}
\definecolor{text1}{cmyk}{1,.65,0,0} 
\definecolor{text2}{rgb}{1,0,0} 
\definecolor{text3}{cmyk}{0,0,0,1} 
\definecolor{text4}{cmyk}{0,0,0,0.5} 

\usepackage{amsfonts}
\usepackage{amsthm}
\usepackage{graphicx}
\usepackage{latexsym}
\makeatletter
\renewcommand{\@seccntformat}[1]
{\csname the#1\endcsname.\enspace}
\makeatother
\newtheorem{theorem}{Theorem}
\newtheorem{lemma}{Lemma}
\newtheorem{remark}{Remark}
\newtheorem{corollary}{Corollary}
\newtheorem{example}{Example}

\begin{document}

\begin{center}
\textbf{On Predictive Density Estimation under $\alpha$-divergence Loss \footnote{ \today}
} 
\\
\end{center}

\begin{center}
{\sc Aziz L'Moudden$^{a}$, \'Eric Marchand$^{a}$, }\\

{\it a  Universit\'e de
    Sherbrooke, D\'epartement de math\'ematiques, Sherbrooke Qc,
    CANADA, J1K 2R1 \quad (e-mails: aziz.lmoudden@usherbrooke.ca; eric.marchand@usherbrooke.ca) } \\

\end{center}
\vspace*{0.2cm}
\begin{center}
{\sc Summary} \\
\end{center}
\vspace*{0.1cm}
\small
Based on $X \sim N_d(\theta, \sigma^2_X I_d)$, we study the efficiency of predictive densities under $\alpha-$divergence loss $L_{\alpha}$ for estimating the density of $Y \sim N_d(\theta, \sigma^2_Y I_d)$.  We identify a large number of cases where improvement on a plug-in density are obtainable by expanding the variance, thus extending earlier findings applicable to Kullback-Leibler loss.  
The results and proofs are unified with respect to the dimension $d$, the variances $\sigma^2_X$ and $\sigma^2_Y$, the choice of loss $L_{\alpha}$; $\alpha \in (-1,1)$.  The findings also apply to a large number of plug-in densities, as well as for restricted parameter spaces with $\theta \in \Theta \subset \mathbb{R}^d$.   The theoretical findings are accompanied by various observations, illustrations, and implications dealing for instance with robustness with respect to the model variances and simultaneous dominance with respect to the loss.

\vspace*{0.5cm}


\section{Introduction}

Consider normally and independently distributed $X|\theta \sim N_d(\theta, \sigma^2_X I_d)$ and $Y|\theta \sim N_d(\theta, \sigma^2_Y I_d)$ and the objective of predicting $Y$ having observed $X$. For predictive analysis purposes, researchers are interested in finding a predictive density $\hat{q}(\cdot \,; X)$ of the density $q(\cdot|\theta)$ of $Y$.  In turn, such a density may play a surrogate role for generating either future or missing values of $Y$. 
To evaluate the performance of such predictive densities, attractive choices are given by the family of $\alpha$-divergence losses (e.g., Csisz\`ar, 1967)
\begin{equation}\label{loss}
L_{\alpha}\left(\theta, \hat{q}\left(\cdot \left| x\right. \right)\right)
=
\int_{\mathbb{R}^d} h_{\alpha}\left(\frac{\hat{q}\left(y \left| x\right. \right)}
{q\left(y \left| \theta \right. \right)}\right) \ q\left(y \left| \theta \right. \right) \,dy,
\end{equation}
where  
\begin{equation}\label{halpha}
h_{\alpha}(z)=\left\{ \begin{array}{ccc}
\frac{4}{1-\alpha^{2}}\left(1-z^{\frac{1+\alpha}{2}}\right) &  & \left|\alpha\right|<1 \\
z \log z &  & \alpha=1 \\
-\log z &   & \alpha=-1\,, \end{array} 
   \right.
\end{equation}
and the corresponding frequentist risk of  $\hat{q}(\cdot;X)$ given by
\begin{equation}\label{risk}
R_{\alpha}(\theta, \hat{q})= E^{X}\left[L_{\alpha}\left(\theta, \hat{q}\left(\cdot \left| X \right. \right)\right) \right]. 
\end{equation}
Notable examples of $L_{\alpha}$ include Kullback-Leibler ($L_{-1}$), reverse Kullback-Leibler ($L_1$), and Hellinger ($L_0$/4).  We point out that it is equivalent to consider 
\begin{equation}\label{halpha-alternative}
h_{\alpha}(z)=\left\{ \begin{array}{ccc}
\frac{4}{1-\alpha^{2}}\left(\frac{1+\alpha}{2} z -z^{\frac{1+\alpha}{2}} + \frac{1+\alpha}{2}\right) &  & \left|\alpha\right|<1 \\
1- z+ z \log z &  & \alpha=1 \\
z-\log z -1 &   & \alpha=-1\,. \end{array} 
   \right.
\end{equation}
With these functions being non-negative, equal to $0$ iff $z=1$,  decreasing for $z \in (0,1)$, increasing for $z>1$, they connect in a more straightforward manner to desirable features for losses in  (\ref{loss}).   The cases $|\alpha| <1$ stand apart, and merit study, as they typically lead to finite loss, unlike the cases of Kullback-Leibler and reverse Kullback-Leibler losses.

This paper is concerned with improvements on plug-in predictive densities of the form
\begin{equation}
\label{plugin}
q_{\hat{\theta}, 1} \sim N_d(\hat{\theta}(X), \sigma^2_Y I_d)\,,
\end{equation}
where $\hat{\theta}(X)$ is a non-degenerate point estimator of $\theta$, and where 
$\theta \in C \subset \mathbb{R}^d$.  Whereas, such plug-in densities were shown in Fourdrinier et al. (2011) to be universally deficient for Kullback-Leibler risk, and improved upon by a subclass of scale expansion variants 
\begin{equation}
\label{pluginc}
q_{\hat{\theta}, c} \sim N_d(\hat{\theta}(X), c^2 \sigma^2_Y I_d)\,,
\end{equation}
with $1<c<c_0$, whereas {\bf all} Bayesian predictive densities are plug-in densities for reverse Kullback-Leibler loss (Yanigomoto and Ohnishi, 2009; Marchand and Sadeghkhani, 2017), analogous analytical results are lacking for the cases $|\alpha| <1$.   This manuscript fills this gap.  Improvements are obtained, interpreted, and illustrated.  The main result (Theorem \ref{main}) reproduces for a large class of choices of $\hat{\theta}$ the Kullback-Leibler result with predictive density estimators $q_{\hat{\theta}, c}$ dominating the plug-in $q_{\hat{\theta}, 1}$ for expansions $c \in (1, c_0)$.   The cut-off point $c_0$ depends on $\alpha, \hat{\theta}$, the parameter space $C$,  the dimension $d$, $\sigma^2_X,$, and $\sigma^2_Y$, but the phenomenon of improvement by variance expansion is otherwise quite general, subject to conditions, namely on $\hat{\theta}$ which we address.  
Further inferences become available as well.  On one hand, we deduce simultaneous dominance results with representatives  $q_{\hat{\theta}, c}$ dominating a given plug-in density $q_{\hat{\theta}, 1}$ for many $\alpha-$ divergence losses, including Kullback-Leibler.  On the other hand, we expand on the fact that dominating predictive densities are not necessarily restricted to normal densities in (\ref{pluginc}), but also include variance mixture of normals with the variance mixing variable taking values on $(1,c_0)$.   The scope of our results is also 
enlarged in view of applications to restricted parameter spaces $C$, including  univariate means constrained to an interval or a half-interval, $L^2$ balls $C=\{\theta \in \mathbb{R}^d: \|\theta\| \leq m\}$, order constraints such as $C=\{\theta \in \mathbb{R}^d:\theta_{i+1} \geq \theta_{i} \hbox{ for } i=1, \ldots d-1 \} $ among many other types (see Remark \ref{constraints}). 

The predictive density estimation framework considered here was put forth for Kullback-Leibler divergence loss in the pioneering work of Aitchison and Dunsmore (1975), as well as Aitchison (1975), and has found applications in information theory, econometrics, machine learning, image processing, and mathematical finance, among others.   
For multivariate normal observables, as considered here, much interest was generated following Komaki (2001) where Bayesian improvements on the minimum risk equivariant predictive density; corresponding to $\hat{\theta}(X)=X$, $c^2=1 + \sigma^2_X/\sigma^2_Y$ in (\ref{pluginc}); were obtained for $d \geq 3$ and Kullback-Leibler loss.   For  Kullback-Leibler loss, further findings with respect to minimaxity, admissibility, dominance, parametric space restrictions, were obtained by George, Liang and Xu (2006), 
Brown, George and Xu (2008), Fourdrinier et al. (2011), among others, while findings for $\alpha-$divergence include the work of Ghosh, Mergel and Datta (2008), Maruyama and Strawderman (2012), Maruyama and Ohnishi (2017), and Marchand and Sadeghkhani (2017). 
The theme of improvements on plug-in densities by variance expansion arises as well for Gamma models (LMoudden et al., 2017) under Kullback-Leibler divergence loss, as well for spherically symmetric and normal models under integrated $L^2$ and $L^1$ losses (Kubokawa, Marchand and Strawderman; 2017, 2015A).  Finally, several researchers have studied the asymptotic efficiency of predictive densities, with $\alpha-$divergence findings for exponential families obtained by Corcuera and Giummol\`{e} (1999).  

The rest of the paper is organized as follows.  At the outset of Section 2, we consider the instructive case $\hat{\theta}(X)=X$ with properties that resonate throughout the manuscript.   In Subsections 2.1 and 2.2, we first consider two specific and instructive cases: {\bf (i)} the affine linear case $\hat{\theta}_a(X)=aX$, and {\bf (ii)} the case 
of a univariate non-negative normal mean and the maximum likelihood estimator $\hat{\theta}_{mle}(X)=\max\{0,X\}$.
In both cases, we obtain necessary and sufficient conditions for a variance expansion $q_{\hat{\theta}, c}$ to dominate the plug-in density $q_{\hat{\theta}, 1}$ and provide numerical illustrations.  We proceed in Subsection 2.3 with a more general dominance finding with implications.  The result is applicable to a large class of plugged-in estimators $\hat{\theta}$ in (\ref{pluginc}), and represents an otherwise unified finding with respect to the dimension $d$, the parameter space $C$, the loss $L_{\alpha}$ for $\alpha \in (-1,1)$, the variances $\sigma^2_X$ and $\sigma^2_Y$.  In Subsection 2.4, we expand on with further analysis and observations relative to the allowable degree of expansion maintaining dominance, as well as simultaneous dominance for several choices of loss $L_{\alpha}$.  We conclude with an example in Subsection 2.5.

\section{Main results}

To begin with, it is instructive to review the case $\hat{\theta}(X)=X$ for the densities in (\ref{pluginc}), and for which the $\alpha-$divergence risk is constant as a function of $\theta \in \mathbb{R}^d$. and given by

\begin{equation}
\label{riska=1}
R_{\alpha}(\theta, q_{\hat{\theta}_1, c}) = \frac{4}{1- \alpha^2} \, \left(1 \, - \, \left(\frac{4 \, c^{1-\alpha}}{2 c^2 (1-\alpha) + 2 (1 + \alpha) + (1-\alpha^2) r}  \right)^{d/2} \,    \right)\,,
\end{equation}
with $r=\sigma^2_X/\sigma^2_Y$ (see Ghosh, Mergel and Datta; 2008; or equation (\ref{riskac}) below that applies for $\hat{\theta}_a(X)=aX$ and $a=1$).  With frequentist risk constant as a function of $\theta$, there exists an optimal choice $c^2_{opt} = 1 + \frac{r (1-\alpha)}{2}$ which minimizes (\ref{riska=1}) in $c^2$.  Observe that this optimal degree of expansion increases in $r$, and decreases in $\alpha$ ranging from a Kullback-Leibler ($\alpha \to -1$) expansion of $1+r$ to an absence of expansion for reverse Kullback-Leibler ($\alpha \to 1$).   Moreover, the risk in (\ref{riska=1}) is for fixed $(d,\alpha,r)$ decreasing in $c$ for $c \in (1,c_{opt})$, and increasing in $c$ for $c>c_{opt}$ (see Theorem \ref{ax}).  This implies that the optimal choice $c^2_{opt}= 1+ \frac{r(1-\alpha_0)}{2}$ for loss $L_{\alpha_0}$ also leads to dominance for Kullback-Leibler loss, as well as losses $L_{\alpha}$ with $\alpha \in (-1,\alpha_0)$.  We thus also have simultaneous dominance with respect to a class of loss functions.  Such features recur throughout the paper in the study of predictive densities (\ref{pluginc}).    We point out that $q_{\hat{\theta}_1, c_{opt}}$ is the (generalized) Bayes predictive density with respect to the prior $\pi(\theta)=1$, as well as the minimum risk equivariant predictive density with respect to changes of location (e.g., Ghosh, Mergel and Datta, 2008).

The amount gained by expanding the variance to the optimal level is reflected by the ratio

\begin{equation}
\label{ratio}
\frac{R_{\alpha}(\theta, q_{\hat{\theta}_1, 1})}{R_{\alpha}(\theta, q_{\hat{\theta}_1, c_{opt}})}\, = \frac{1 - \left\lbrace1+ (r(1-\alpha^2)/4)\right\rbrace^{-d/2}}{\;\;\;\;\;\;\; 1- \left\lbrace(1+ (r(1-\alpha)/2)) \right\rbrace^{-d(1+\alpha)/4}}\,.
\end{equation}

Numerical and analytical evaluations (see Figure 1) suggest that this ratio decreases in $\alpha \in (-1,1)$, as well as in $d$.  As a function of $r$, the ratio approaches $1$ as $r \to 0$ and $r \to \infty$, increasing for $r < r_0$ up to a maximal $r=r_0$, and decreasing for $r>r_0$.  As an exemplar, for $\alpha=0, d=2$, the ratio above reduces to $\frac{2+r + \sqrt{4+2r}}{4+r}$, and behaves as above with a maximum value of around $1.2071$ attained at $r_0 = 4(1+\sqrt{2}) \approx 9.6568$.  Finally, further numerical evaluations  suggest that such maximal gains increase as $\alpha$ decreases and are attenuated with increasing dimension $d$.  

As mentioned above, it seems plausible that such aspects of this benchmark case may recur for other choices $\hat{\theta}$ in predictive densities (\ref{pluginc}).

\begin{figure}[!ht]
  \centering
    \includegraphics[width=0.99\textwidth]{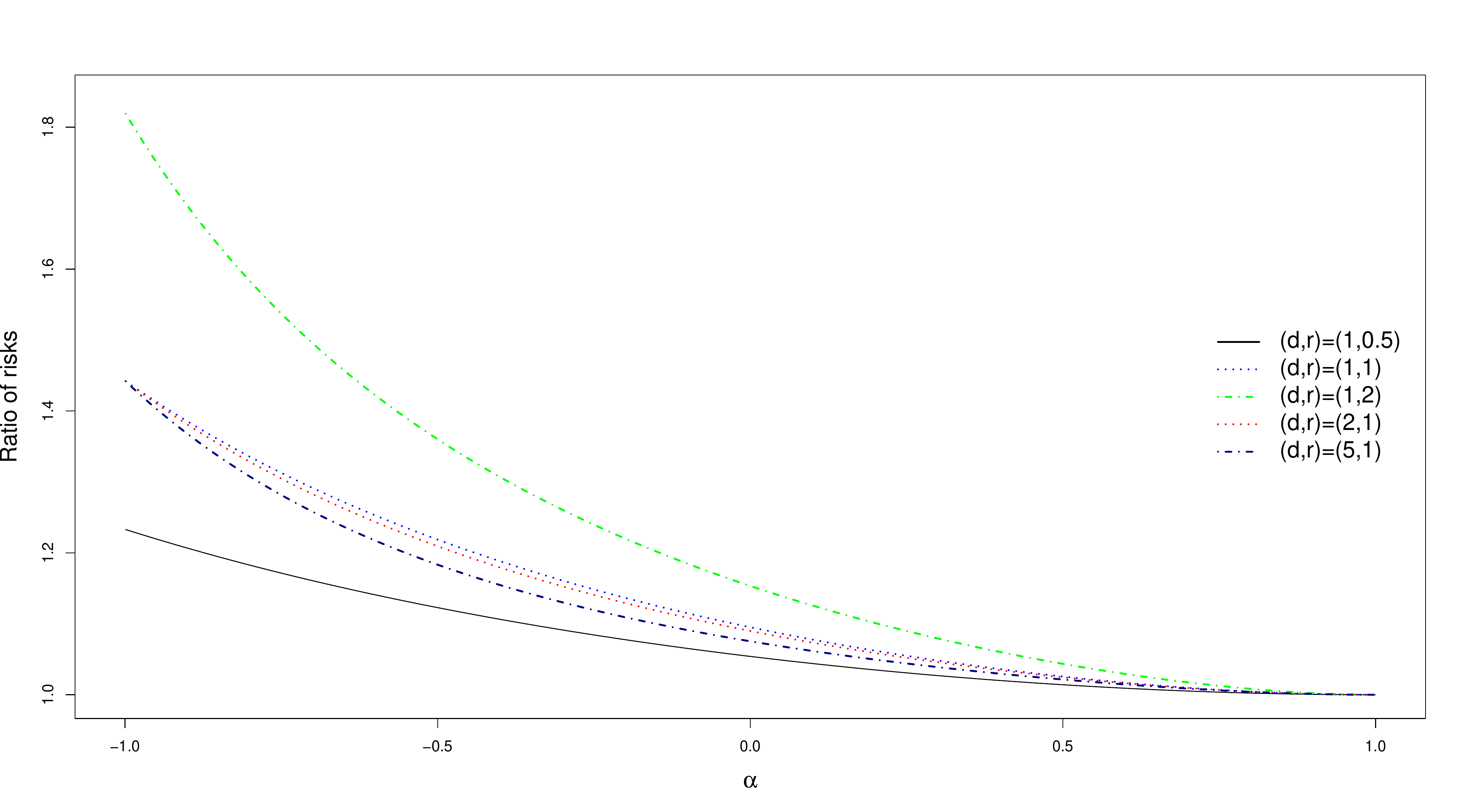}
     \label{a=1}
  \caption{Risk ratios (\ref{ratio}) as functions of $\alpha$ for various pairs $(d,r)$.
  }
\end{figure}

\subsection{Case $\hat{\theta}_a(X)=a X$}

Before analyzing the risk performance of predictive densities in (\ref{pluginc}) associated with affine linear estimators, we will first require the following result, which highlights the relationship between $\alpha-$divergence loss and reflected normal loss.  The result is known (e.g., Ghosh, Mergel and Datta; 2008) and we will expand on its significance below in Remark \ref{significance}.  Hereafter, we denote $\phi$ as the $N_d(0,I_d)$ p.d.f., and $\Phi$ as the $N(0,1)$ c.d.f.

\begin{lemma}  
\label{alphadual}   
For $|\alpha|<1$, the loss $L_{\alpha}$ incurred by the predictive density estimate $q_{\hat{\theta},c}$ as in (\ref{pluginc}) for estimating the density $q(\cdot|\theta)$ of $Y$ is given by:

\begin{equation}
\label{L(q,c)}
L_{\alpha}(\theta,q_{\hat{\theta},c})=\frac{4}{1-\alpha^2}\left( 1 - \left(\dfrac{2\,\,c^{1-\alpha}}{(1-\alpha)\,c^2+1+\alpha} \right)^{\tfrac{d}{2}}exp\left( -\dfrac{||\hat{\theta}(x)-\theta||^2}{2\gamma_0} \right)  \right),  
\end{equation} 
with $\gamma_0 = 2 \left(  \frac{c^2}{1+\alpha} + \frac{1}{1-\alpha} \right) \sigma_Y^2 .$

\end{lemma}
{\bf Proof.}  The result follows by a development of (\ref{loss}) for 
$\hat{q}(y;x)\,=\, \frac{1}{(c \sigma_Y)^d} \phi(\frac{y - \hat{\theta}(x)}{c \sigma_Y})$ and $q(y|\theta) \, = \,  \frac{1}{(c\sigma_Y)^d} \phi(\frac{y - \theta}{\sigma_Y})$, $y \in \mathbb{R}^d$.  \qed 

\begin{theorem}
\label{ax}
Consider $X\sim N_d(\theta, \sigma^2_X I_d)$ independent of  $Y\sim N_d(\theta, \sigma_Y^2 I_d)$ and the problem of estimating the density of $Y$ under $\alpha-$divergence loss $L_{\alpha}$ as in (\ref{loss}) with $|\alpha|<1$ and $\theta \in \mathbb{R}^d$.   Let $\hat{\theta}_a(X)=aX$ for $0 < a < 1$, set $r=\sigma^2_X/\sigma^2_Y$, and consider predictive density estimators $q_{\hat{\theta}_a, c}(\cdot;X) \sim N_d( aX, c^2 \sigma^2_Y I_d)$ for $c \geq 1$.  Then $q_{\hat{\theta}_a, c}$ dominates the plug-in density $q_{\hat{\theta}_a, 1}$ if and only if $1 < c^2 \leq k^2(\alpha, a, r)$, where $k(\alpha, a, r)$ is the solution in $c \in (1,\infty)$ of 
\begin{equation}
\label{K}  
2\, (1-\alpha)\, c^2 -\left(4+(1-\alpha^2)\,a^2\,r \right) \, c^{1-\alpha}+(1+\alpha) \, \left(2+(1-\alpha)\,a^2\,r \right) \, = \, 0 \,.
\end{equation}
\end{theorem}
{\bf Proof.}  It follows from Lemma \ref{alphadual} that 
\begin{equation}
\nonumber
R_{\alpha}(\theta, q_{\hat{\theta}_a, c}) = \frac{4}{1- \alpha^2} \left( 1 \, - \, 
\left(\frac{4 \sigma^2_Y c^{1-\alpha}}{\gamma_0 (1-\alpha^2)} \right)^{d/2} \, \mathbb{E} 
\left(e^{- \frac{\|aX-\theta\|^2}{2 \gamma_0}} \right) \right)\,, 
\end{equation}
with $\gamma_0= 2 \sigma^2_Y \frac{c^2 (1-\alpha) + 1+ \alpha }{1-\alpha^2}\,. $
Since $\|aX-\theta\|^2 \sim a^2 \sigma^2_X \, \chi^2_d (\frac{(a-1)^2 \|\theta\|^2}{a^2 \sigma^2_X})\,$, one obtains
\begin{equation}
\nonumber
\mathbb{E}(e^{- \frac{\|aX-\theta\|^2 }{2 \gamma_0}}) \, = \,  (\frac{\gamma_0}{\gamma_0 + a^2 \sigma^2_X})^{d/2} \, e^{- \frac{(a-1)^2 \|\theta\|^2}{2 (\gamma_0 + a^2 \sigma^2_X)}}\,,
\end{equation}
whence the expression
\begin{equation}
\label{riskac}
R_{\alpha}(\theta, q_{\hat{\theta}_a, c}) = \frac{4}{1- \alpha^2} \, \left(1 \, - \, \left(\frac{4 \, c^{1-\alpha}}{2 c^2 (1-\alpha) + 2 (1 + \alpha) + a^2 (1-\alpha^2) r}  \right)^{d/2} \, e^{\, - \frac{(a-1)^2 \|\theta\|^2}{2 (\gamma_0 + a^2 \sigma^2_X)}}\,   \right)\,.
\end{equation}
Now, observe for the difference in risks that, for $c>1$,
\begin{eqnarray}
\nonumber & & \frac{1-\alpha^2}{4} \{ R_{\alpha}(\theta, q_{\hat{\theta}_a, c}) - R_{\alpha}(\theta, q_{\hat{\theta}_a, 1}) \} 
\\
\label{2.9} & \leq & \left\lbrace \left(\frac{4}{4+a^2 (1-\alpha^2) r}    \right)^{d/2} \, - \, \left( H(c)\right)^{d/2} \right\rbrace \; e^{\, - \, \frac{(a-1)^2 \|\theta\|^2/ 2\sigma^2_Y}{a^2 r + 4/(1-\alpha^2)} }\,,
\end{eqnarray}
with 
\begin{equation}
\label{H}
H(c) \, = \, \left(\frac{4 \, c^{1-\alpha}}{2 c^2 (1-\alpha) + 2 (1 + \alpha) + a^2 (1-\alpha^2) r}    \right)\,,
\end{equation}
and with equality if and only if $\theta=0$.  The result follows by verifying that, for $c>1$, $H(c) \geq \frac{4}{4+a^2 (1-\alpha^2) r}$ if and only if  to $1 < c^2 \leq k^2(\alpha, a, r)$.  \qed

\begin{remark}
\label{significance}
This paper focuses on the effect of variance expansion, that is the role of $c^2$ on the frequentist risk performance of $q_{\hat{\theta},c}$.   Alternatively, it is natural and of interest to study the role of the plugged-in estimator $\hat{\theta}$.   In view of expression (\ref{L(q,c)}), it is apparent that the frequentist risk under $\alpha-$divergence loss of the predictive density estimator $q_{\hat{\theta},c}$ relates to the point estimation risk performance of $\hat{\theta}$ as an estimator of $\theta$ under reflected normal loss $L(\theta, \hat{\theta}) \, = \, 1 - e^{- \|\hat{\theta}-\theta\|^2/2\gamma_0}$.  Ghosh, Mergel and Datta (2008) capitalized on such a dual relationship to derive predictive densities $q_{\hat{\theta},c}$ dominating the minimum risk equivariant, and minimax, predictive density $\hat{q}_{mre}(\cdot; X) \sim N_d(X, (\frac{1-\alpha}{2} \sigma^2_X +\sigma^2_Y) I_d)$ for $d \geq 3$  Further applications were recently obtained by Marchand, Perron and Yadegari (2017) for $\alpha-$divergence prediction, while additional point estimation results for reflected normal loss were obtained by Kubokawa, Marchand and Strawderman (2015A, 2015B). 
\end{remark}

Here are some further observations and implications of Theorem \ref{ax}.

\begin{remark}
\label{remarkTheorem2.1}
\begin{itemize}

\item  For Hellinger loss (i.e., $\alpha=0$), Theorem \ref{ax}'s cut-off point simplifies to $k^2(0,a, r)= (1+\frac{a^2r}{2})^2.$   In this case, and more generally for other choices of $\alpha$, it is easy to show that $k(\alpha, a, r)$ increases in $r=\sigma^2_X/\sigma^2_Y$, converging to $1$ as $r \to 0$.  Smaller values of $r$; which may also translate to larger samples sizes from $X$;  correlate with greater efficiency of $aX$ for estimating $\theta$ and less of a need to expand on the plug-in density.  Larger values of $r$ have the opposite effect.  One can also infer a robustness result :  if the ratio $r$ of variances is misspecified and that the actual ratio is equal to $r'$ dominance of $q_{\theta_a,c}$ over $q_{\theta_a,c}$ persists for $1 < c^2 \leq k^2(\alpha, a, r)$ as long as $r'>r$, i.e., one has underestimated the ratio of variances.  

\item  An explicit lower bound for Theorem \ref{ax}'s cut-off point is given by the inequality $k^2(\alpha, a, r) \geq 1 + \frac{a^2r (1-\alpha)}{2}$.  This follows from (\ref{2.9}) and observing that $H(c)$ in (\ref{H}) increases in $c$ for $1 < c^2 \leq 
1 + \frac{a^2r (1-\alpha)}{2}$.  

\item  Theorem \ref{ax} applies for $a=1$, with the $1 <c^2 < k^2(\alpha,1,r)$ as the necessary and sufficient condition for dominance.

\item  Numerical evidence suggests that, for fixed $a \in (0,1)$, $r>0$,  $k(\alpha, a, r)$ decreases in $\alpha$, which is quite plausible (and undoubtedly true for the lower bound in the previous paragraph).  If true, choices $c^2=k^2(\alpha_0, a, r)$ would not only lead to domination for loss $L_{\alpha_0}$, but also for all other $\alpha-$divergence losses with $-1 \leq \alpha \leq \alpha_0$.

\end{itemize}
\end{remark}

We conclude this section with a numerical illustration.

Figure 2 represents for Hellinger loss (i.e., $\alpha=0$), $d=3$, $\sigma^2_X=1$, $a=0.75$, ratios of risks $\frac{R_0(\theta, q_{\hat{\theta}_0.75,c})}{R_0(\theta, q_{\hat{\theta}_0.75,1})}$, as a function of $\|\theta\|$, for $r=0.5,1,2$ and $c=k(\alpha,a,r)$ and $c=(1+k(\alpha,a,r))/2$.  The graphs illustrate the dominance result given by Theorem \ref{ax} and permit us to focus here on the effect of the ratio of variances $r=\frac{\sigma^2_X}{\sigma^2_Y}$, as well as the degree of variance expansion $c^2$, in relationship to the gains that are attainable by variance expansion as opposed to the plug-in density.   Overall, the gains can be significant, as illustrated here for a specific setting of $\alpha, d, \sigma^2_X$, and $a$.   Here are some observations based on  Figure 2 and several other numerical evaluations. 

\begin{enumerate}
\item[ {\bf (i)}]  For larger $r$, maximal gains are more important, as well as gains for small or moderate $\|\theta\|$.  This is a recurrent feature below in other situations, for other plug-in choices $\hat{\theta}(X)$.  Given the decreasing relative reliability of the information provided by $X$ for making inferences about $\theta$, such behaviour is somewhat anticipated and relates in this example to the allowable degree of expansion to maintain dominance which increases in $r$ (Remark \ref{remarkTheorem2.1}).  It is somewhat delicate here as the choice itself of $a$ would typically depend on $\sigma^2_X$, as is the case of the posterior expectation for a $N_d(0,I_d)$ prior corresponding to the multiple $\frac{d}{d+\sigma^2_X} $ ; 

\item[ {\bf (ii)}]   Expansions to $c=k(\alpha,a,r)$ lead to more important maximal gains, while the compromise choice $c=(1+k(\alpha,a,r))/2$ flattens out the gains with better performance for small $\|\theta\|$.  Observe as well the ratios equal $1$ for $\theta=0$ as established within the proof of Theorem \ref{ax}, and that more important expansions will not lead to dominance in view of the necessity and sufficiency of Theorem \ref{ax}.

\end{enumerate}
 
\begin{figure}[!ht]
  \centering
    \includegraphics[width=0.99\textwidth]{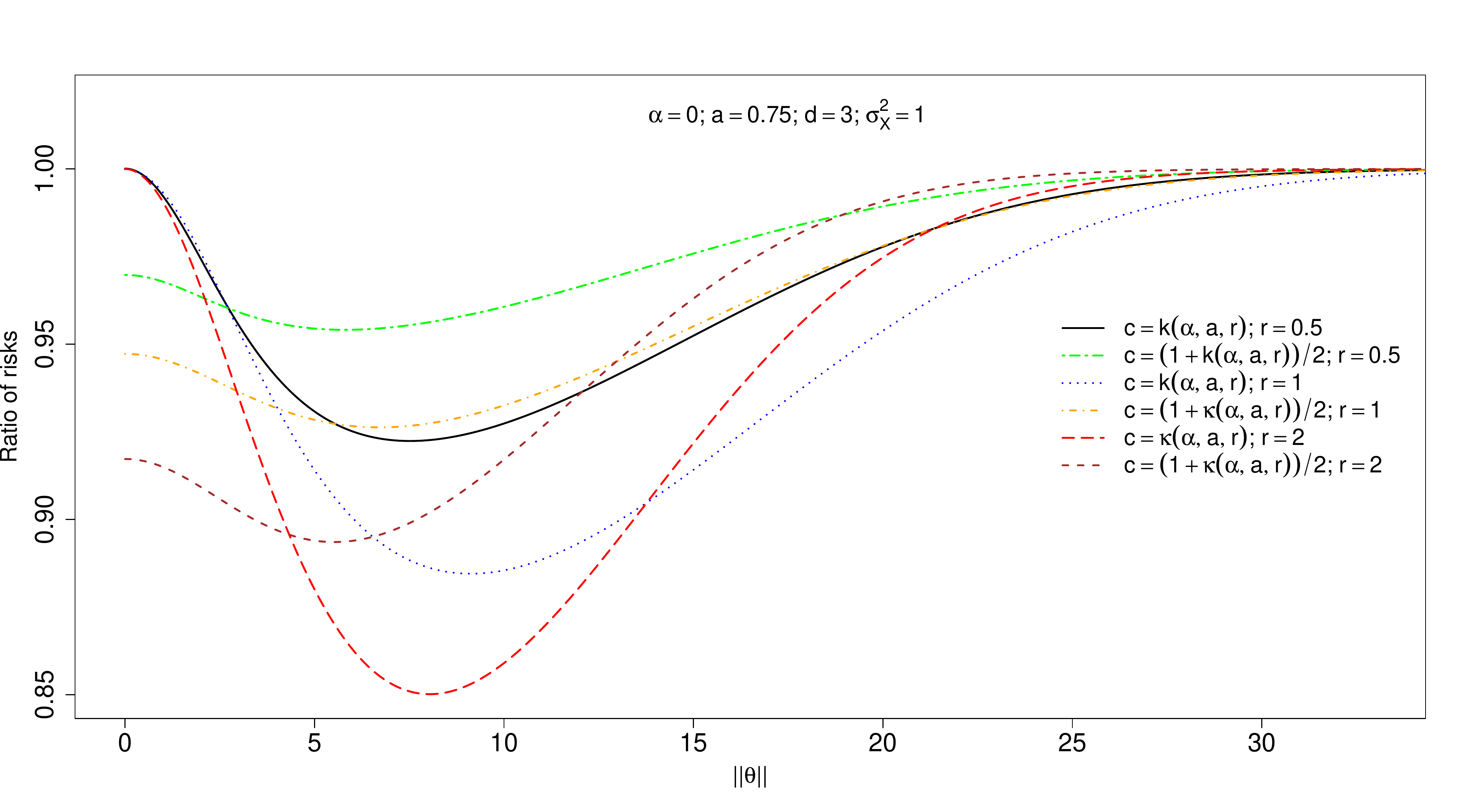}
     \label{ab1}
  \caption{Hellinger risk ratios $\frac{R_0(\theta,q_{\hat{\theta}_{a},c})}{R_0(\theta,q_{\hat{\theta}_{a},1})}$ for various $c$ and $r$}
\end{figure}

\subsection{Case of a non-negative mean with $\hat{\theta}_+(X)=\max\{X,0\}$}

As in the previous section, we proceed with an instructive example bringing into play a non-negativity constraint and the choice of the maximum likelihood estimator (mle) $\hat{\theta}_+(X)=\max\{X,0\}$.  The relative tractability of the $\alpha$-divergence
risk, which arises with a convenient expression for the expected reflected normal loss, leads to a necessary and sufficient condition for a variance expansion to dominate the predictive mle.

\begin{theorem}
\label{x+}
Consider $X\sim N(\theta, \sigma^2_X)$ independent of  $Y\sim N(\theta, \sigma_Y^2)$ and the problem of estimating the density of $Y$ under $\alpha-$divergence loss $L_{\alpha}$ as in (\ref{loss}) with $|\alpha|<1$ and $\theta \geq 0$.   Let $\hat{\theta}_+(X)= \max \{X,0\}$, and consider predictive density estimators $q_{\hat{\theta}_+, c}(\cdot;X) \sim N(\hat{\theta}_+(X) , c^2 \sigma^2_Y)$ for $c \geq 1$. Set $r=\sigma^2_X/\sigma^2_Y$, $\gamma_0(c)= 2 \sigma^2_Y \left(\frac{c^2(1-\alpha) + 1+ \alpha}{1-\alpha^2} \right)$, $\gamma_1(c)= \sqrt{\frac{\gamma_0(c)}{\gamma_0(c)+\sigma^2_X}}$, $A_1(c) = \sqrt{\frac{2 c^{1-\alpha}}{c^2(1-\alpha) + 1+ \alpha}}$, and $A_2(c)=-1+A_1(c) (1+\gamma_1(c))$.
Then $q_{\hat{\theta}_+, c}$ dominates the plug-in m.l.e. density  $q_{\hat{\theta}_+, 1}$ if and only if $1 < c \leq \kappa(\alpha, r)$, where $\kappa(\alpha, r)$ is the solution in $c \in (1,\infty)$ of 
\begin{equation}
\label{kappa}
A_2(c) \, = \, \left(1+ \frac{r (1-\alpha^2)}{4}  \right)^{-1/2}\,.
\end{equation}
The dominance is strict for $\theta>0$ or $1<c<\kappa(\alpha,r)$, with equality of risks if and only of $\theta=0$ and $c=\kappa(\alpha,r)$.
\end{theorem} 
{\bf Proof.} {\bf (I)}  It follows from Lemma \ref{alphadual} that
\begin{equation} 
\label{risktheta+}
R_{\alpha}(\theta, q_{\hat{\theta}_+, c}) \, = \, \frac{4}{1-\alpha^2} \left[1- A_1(c) \, G(\theta,c)  \right]\,, \end{equation}
with  $G(\theta, c) \,=\, \mathbb{E} \left(e^{- \frac{|\hat{\theta}_+(X) - \theta|^2}{2  \gamma_0(c)}}  \right)$ and the given notation.  Calculations yield the expression $G(\theta, c) \, = \, G_1(\theta, c) \, + \, G_2(\theta, c)$, with $G_1(\theta, c)\, = \, e^{- \frac{\theta^2}{2 \gamma_0(c)}} \, \Phi(- \frac{\theta}{\sigma_X})$ and $G_2(\theta, c)\, = \, \gamma_1(c) \, \Phi(\frac{\theta}{\gamma_1(c) \sigma_X}) $, as well as 
$$  \frac{\partial}{\partial \theta} G(\theta, c) \, = \, -\frac{\theta}{\gamma_0(c)} \, G_1(\theta, c)\,.$$    

{\bf (II)} For the difference in risks, we thus obtain from the above 
\begin{eqnarray}  
\nonumber \Delta(\theta,c) \, & = & \, \frac{1-\alpha^2}{4} \left\lbrace R_{\alpha}(\theta, q_{\hat{\theta}_+, c}) \, - \, R_{\alpha}(\theta, q_{\hat{\theta}_+, 1}) \right\rbrace \\
\label{Delta+}
\, & = & G(\theta, 1) \, - \, A_1(c) \, G(\theta, c)\,,
\end{eqnarray}
and
\begin{equation}
\label{derivativeofdelta}
\frac{\partial}{\partial \theta} \, \Delta(\theta, c) \, = \,
\frac{\theta  e^{- \frac{\theta^2}{2 \gamma_0(1)}} \Phi(-\frac{\theta}{\sigma_X})}{\gamma_0(1)} \left\lbrace  -1 \; + \; A_1(c) \frac{\gamma_0(1)}{\gamma_0(c)} \, e^{\frac{\theta^2}{2}(\frac{1}{\gamma_0(1)} - \frac{1}{\gamma_0(c)} )} \right\rbrace  \,.
\end{equation}
Now, since $\gamma_0(c)$ increases in $c$, $c \geq 1$, and $A_1(c) \leq 1$ for all $c \geq 1$ with equality iff $c=1$, we infer, for fixed $c \in (1,\infty)$, that $\frac{\partial}{\partial \theta} \, \Delta(\theta, c)$ changes signs from $-$ to $+$ as $\theta$ varies on $[0,\infty)$.  Therefore, for a given $c >1$, $q_{\hat{\theta}_+, c}$ will dominate $q_{\hat{\theta}_+, 1}$ if and only if $\Delta(0,c) \leq 0$ and $\lim_{\theta \to \infty} \Delta(\theta,c) \leq 0$.  Furthermore, from (\ref{Delta+}) and the earlier expression for $G(\theta,c)$, we obtain
$$ \Delta(0,c) - \lim_{\theta \to \infty} \Delta(\theta,c) \, = \frac{1- \gamma_1(1)}{2} - A_1(c) (\frac{1- \gamma_1(c)}{2}) \, \geq 0 \,, $$
since $A_1(c) \leq 1$ and $ 1 > \gamma_1(c) \geq \gamma_1(1)$ for all $c \geq 1$. \\

{\bf (III)}
We thus have that $q_{\hat{\theta}_+, c}$ will dominate $q_{\hat{\theta}_+, 1}$ if and only if $\Delta(0,c) \leq 0$, and there remains to show that this inequality is equivalent to the stated condition, and with equality if and only if $c=\kappa(\alpha,r)$.  To justify this last step, 
since $\lim_{c \to \infty} \Delta(0,c) \,=\, G(0,1)>0$ from (\ref{Delta+}),
it suffices to show that  $\Delta(0,c)$ decreases and then increases, as a function of $c \in 
(1,\infty)$.

In turn, by virtue of (\ref{Delta+}), it will suffice to show that 
\begin{equation}
sgn \left( \frac{\partial}{\partial c} A_1(c) G(0,c) \right) \, \hbox{ varies from  + to - for } c \in (1,\infty).
\end{equation}
Setting $H(c) = \sqrt{c^2(1-\alpha) + (1+\alpha)}$, we may express
$$  A_1(c) G(0,c) \, = \, \frac{1}{\sqrt{2}} \; c^{\frac{1-\alpha}{2}} \left(\frac{1}{H(c)} + \frac{1}{\sqrt{H^2(c)+ \frac{r (1-\alpha^2)}{2}}}  \right)\,.$$
With $H'(c) \, = \, \frac{(1-\alpha) c}{H(c)}$, some calculations and manipulations permit us to write
\begin{equation}
\nonumber
\frac{\partial}{\partial c} A_1(c) G(0,c) \, = \, - \frac{(1-\alpha^2) c^{-\frac{1+\alpha}{2}}}{(\sqrt{2}H(c))^3} \; \left(c^2-1 + W(c)\right)\,,
\end{equation}
with $$W(c) \, = \, (\frac{c^2-1- r(1-\alpha)/2}{(H^2(c) + r(1-\alpha^2)/2)^{3/2}}) \; H^3(c)\,.$$  From this, we see that $\lim_{c \to 1^+} \, (c^2-1+W(c))<0$, implying that $\frac{\partial}{\partial c} A_1(c) G(0,c)$ is positive for small enough $c$.  On the other hand, since $c^2-1+W(c) \geq 0$ for $c^2 \geq 1+ r(1-\alpha)/2$, it will suffice to complete the proof that $W(c)$ be increasing in $c$ for $c^2 < 1+ r(1-\alpha)/2$.  Finally, a calculation yields the expression
\begin{eqnarray*} \frac{\partial}{\partial c} \, \log W(c) \, &=&\, 
\frac{\partial}{\partial c}  \left\lbrace
3 \log B(c)  + \log  \left(\frac{r(1-\alpha)}{2} + 1 - c^2 \right) - 
\frac{3}{2} \log \left( B^2(c) + \frac{r(1-\alpha^2}{2} \right) \right\rbrace \\
& \, = &  \, 3c(1-\alpha) + \frac{2c}{\frac{r(1-\alpha)}{2} + 1 - c^2 } -  \frac{3c(1-\alpha)/2}{H^2(c) + \frac{r(1-\alpha^2)}{2}}\, \\
& \geq & 0 \,,
\end{eqnarray*}
for $c^2 < 1+ r(1-\alpha)/2$, which establishes the result.   \qed

\begin{remark}
\label{remarkmax(0,x)}
As a function of $\theta$, $\theta \geq 0$, the frequentist risks $R_{\alpha}(\theta,q_{\hat{\theta}_+,c})$ are increasing with a limiting value at $\theta \to \infty$ equal to $\frac{4}{1-\alpha^2} \left( 1-A_1(c) \gamma_1(c) \right)$, and with a value of $\frac{4}{1-\alpha^2} \left( 1-\frac{A_1(c)}{2} (1+\gamma_1(c)) \right)$ at $\theta=0$.  These properties are obtained from (\ref{risktheta+}).
As in Remark \ref{remarkTheorem2.1}, the cut-off point $\kappa(\alpha,r)$ can be shown to be increasing as a function of $r$, and is decreasing as a function of $\alpha$ according to numerical evaluations.  The former is obtained in continuity with the arguments of the proof of Theorem \ref{x+} and with the r.h.s. of (\ref{kappa}) decreasing in $r$.  

\end{remark}

\begin{remark}
\label{robustness-negativetheta}
A surprising robustness result is also available from the analysis above in the proof of Theorem \ref{x+}.  Indeed, it is also the case that the difference in risks $\Delta(\theta,c))$ is negative for all $\theta<0$ and $1<c \leq \kappa(\alpha,r)$.  In other words, the expansions $q_{\hat{\theta}_+,c}$ that dominate $q_{\hat{\theta}_+,1}$ with lower risk on $[0,\infty)$ continue providing lower $\alpha-$divergence frequentist risk for negative values of $\theta$.  This is relevant to cases where it is believed that the constraint $\theta \geq 0$ holds true, but, unbeknownst to the investigator, the actual value of $\theta$ is negative.  In such cases, the choice of the plug-in estimate $\max\{x,0\}$ is, of course, not desirable, but the expansion offers better protection against the misspecification.  Finally, to see why the difference is risks remains negative for negative $\theta$; for all $r>0$ and choice of loss $L_{\alpha}$ for $|\alpha|<1$; it suffices to observe from (\ref{derivativeofdelta}) that the sign of $\frac{\partial}{\partial \theta} \Delta(\theta,c)$ varies from $-$ to $+$ as $\theta$ increases from $-\infty$ to $0$ and that $\lim_{\theta \to -\infty} \Delta(\theta,c) = 0$ for all $c \geq 1$ as seen directly by working with (\ref{risktheta+}).  Accordingly, analogous inferences with respect to a persistent dominance result when underestimating $r$, as well as a simultaneous dominance result for various choices of $L_{\alpha}$, apply.
\end{remark}

\begin{example}
We conclude this section with a numerical illustration.   Figure 3 exhibits the relative frequentist risk performance of the maximum likelihood density $q_{\hat{\theta}_+,1}$ and the variance expansion $q_{\hat{\theta}_+,c}$ with Theorem \ref{x+}'s cut-off point $c=\kappa(\alpha,r)$.
More specifically, risk ratios for $\sigma^2_X=1$ are drawn for various combinations of $r$ and $\alpha$.  Theoretically, the ratios are bounded by $1$, and we point out the equality of risks at $\theta=0$ in accordance with Theorem \ref{x+}.  As shown by the graphs, the gains can be significant, tend to be more important for smaller values of $\alpha$, and large values of $r$.   Other levels of variance expansion, such as $c^2=(1+\kappa(\alpha,r))/2$ have the same effect as in Figure 2.
\end{example}
 
\begin{figure}[!ht]
  \centering
    \includegraphics[width=0.99\textwidth]{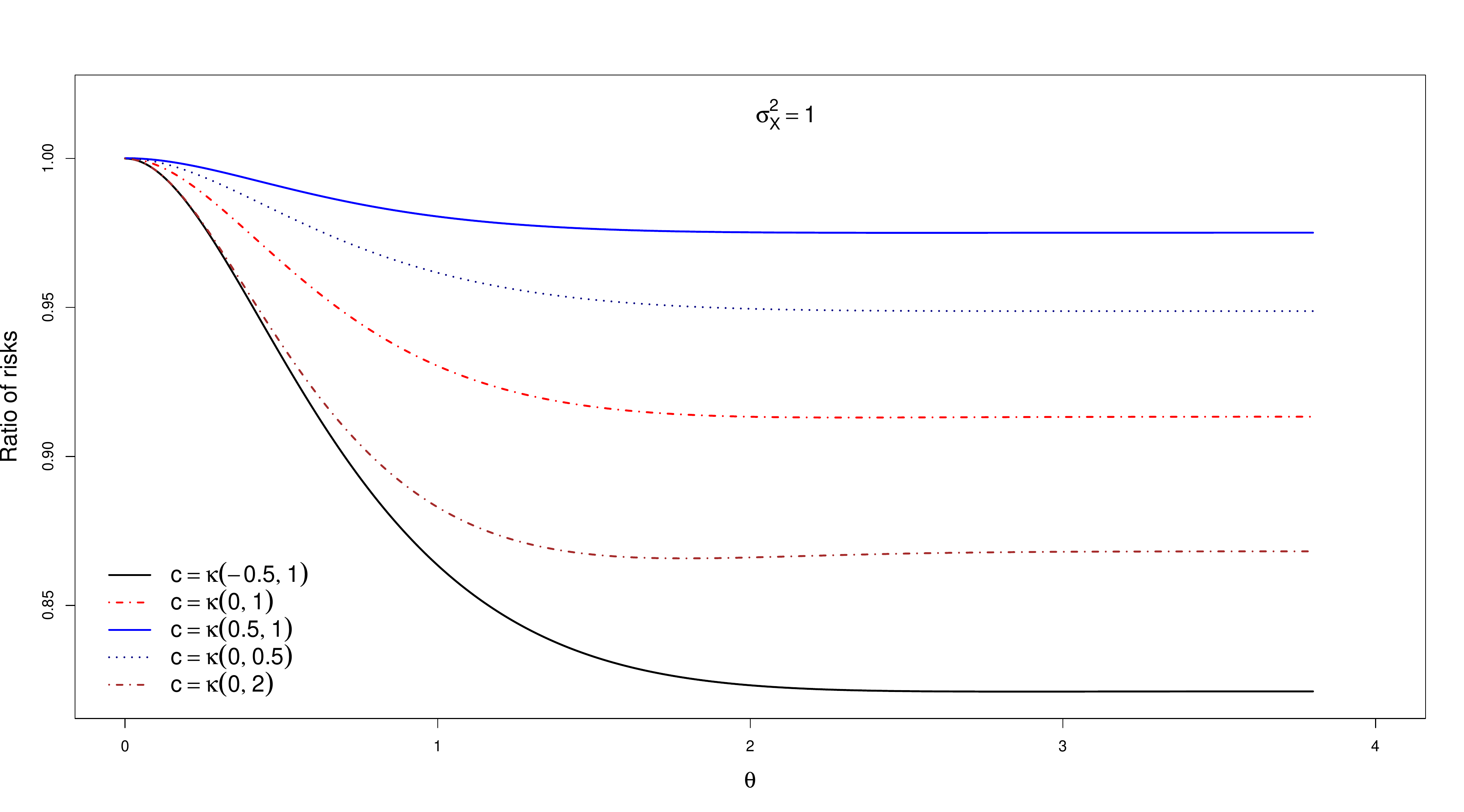}
     \label{ab2}
  \caption{Risk ratios  $\frac{R_{\alpha}(\theta,q_{\hat{\theta}_+,c})}{R_{\alpha}(\theta,q_{\hat{\theta}_+,1})}$ for $c=\kappa(\alpha,r)$, and various $c$ and $r$}
\end{figure}

\subsection{General $\hat{\theta}(X)$}

We begin with the following which we will require.

\begin{lemma}
\label{bill}
Let $T$ be a non-negative and continuous random variable such that  $\mathbb{E}(T^2) < \infty$.   Let $s$ be a positive constant.   Then, we have 
\begin{equation}
\label{billA}
  \mathbb{E}(T e^{-sT}) \geq  \mathbb{E}(T) \, e^{-s \, \mathbb{E}(T^2)/\mathbb{E}(T)} \,.
\end{equation}  
    Moreover, if the distribution of $T$ depends on a parameter $\theta \in C$, and if there exist positive constants $b_0, b_1, b_2$ such that $ b_0 \leq \mathbb{E}_{\theta}(T) \leq b_1$ and $ \mathbb{E}_{\theta}(T^2) \leq b_2 $ for all $\theta \in C$, then we have 
\begin{equation}
\label{billB}
  \mathbb{E}_{\theta}(T e^{-sT}) \geq b_0 \, e^{-s \frac{b_2}{b_1}} \,,
\end{equation}  
for all $\theta \in C$.
\end{lemma}
{\bf Proof.}  Let $f_T$ be the density of $T$ and let $W$ be a random variable with density $w f_T(w)/\mathbb{E}(T)\,, w>0$.  We then have  
$$  \mathbb{E}(T e^{-sT}) \,=\, \mathbb{E}(T) \, \mathbb{E}(e^{-sW})  \geq \mathbb{E}(T) \, (e^{-s \, \mathbb{E}W}) \, = \mathbb{E}(T) \, e^{-s \, \frac{\mathbb{E}(T^2)}{\mathbb{E}(T)}}\,,$$
by using Jensen's inequality.  This establishes (\ref{billA}), and (\ref{billB}) is a direct consequence of (\ref{billA}).       \qed

We now are ready for our main result.

\begin{theorem}
\label{main}
Consider $X\sim N_d(\theta, \sigma^2_X I_d)$ independent of  $Y\sim N_d(\theta, \sigma_Y^2 I_d)$ and the problem of estimating the density of $Y$ under $\alpha-$divergence loss $L_{\alpha}$ as in (\ref{loss}) with $|\alpha|<1$ with $\theta \in C$.  Let $q_{\hat{\theta},1}$ be a plug-in predictive density based on a non-degenerate $\hat{\theta}(X)$, let $ Z=
\frac{\|\hat{\theta}(X)-\theta\|^2 }{\sigma_Y^2}$, and define $\tau= (1-\alpha) \epsilon - 2 \alpha d$,  with  
\begin{equation}
\label{con11}
\epsilon \equiv \epsilon(\alpha)= \inf_{\theta \in C} \mathbb{E}_{\theta} (Z \, e^{- (1-\alpha^2) Z/8})\,.
\end{equation}
\begin{enumerate}
\item[ {\bf (a)}]
Then, assuming $\epsilon>0$, the expansion $q_{\hat{\theta},c} \sim N_d(\hat{\theta}(X), c^2 \sigma^2_Y I_d)$ dominates $q_{\hat{\theta},1}$ whenever 
$$1 < c^2 \leq k(d,\alpha, \sigma^2_X, \sigma^2_Y) =  \frac{\tau + \sqrt{\tau^2 + 4 d^2(1-\alpha^2)}}{2d(1-\alpha)}\,\,.$$
\item[ {\bf (b)}]  Furthermore, if  $b_0, b_1, b_2$ are positive numbers such that 
$b_0 \sigma^2_Y \leq \mathbb{E}_{\theta} \|\hat{\theta}(X)-\theta\|^2 \leq b_1 \sigma^2_Y$ and $\mathbb{E}_{\theta} \|\hat{\theta}(X)-\theta\|^4 \leq b_2 \sigma^4_Y$ for all $\theta \in C$, then the assumption $\epsilon >0$ is satisfied and 
\begin{equation}
k(d,\alpha, \sigma^2_X, \sigma^2_Y) \geq \frac{\underline{\tau} + \sqrt{\underline{\tau}^2 + 4 d^2(1-\alpha^2)}}{2d(1-\alpha)}\,,
\end{equation}
with $\underline{\tau} = (1-\alpha)\, b_0 \,
e^{- \frac{1-\alpha^2}{8} \frac{b_2}{b_1} }\, - 2\alpha d\,$.
\end{enumerate}
\end{theorem}
{\bf Proof.}  For part {\bf (a)}, setting $B(c)=(1-\alpha)c^2+(1+\alpha)$, $H_{\theta}(c)=\mathbb{E}_{\theta}(e^{-\frac{(1-\alpha^2)}{4B(c)}Z})$, and $h_{\theta}(c)=\frac{c^{1-\alpha}}{B(c)} \ H^{\frac{2}{d}}_{\theta}(c)$, we obtain from Lemma \ref{alphadual} the risk expression:

\begin{equation}
\label{risk}
R_{\alpha}(\theta, q_{\hat{\theta},c})=\frac{4}{1-\alpha^{2}} \left(1-2^{\frac{d}{2}} \ h_{\theta}^{\frac{d}{2}}(c)\right).
\end{equation}
To establish the result, it will suffice to show that, for all $\theta \in C$, 
 $\frac{d}{dc}h_{\theta}(c)> 0$ for  $c\in (1,k(d,\alpha, \sigma^2_X, \sigma^2_Y)]$;
i.e., the risk $R_{\alpha}(\theta, q_{\hat{\theta},c})$ decreases, for all $\theta \in C$, as a function of $c$, for $1<c^2< k(d,\alpha, \sigma^2_X, \sigma^2_Y)$. 
We have
  \begin{equation}
  \label{hderivative}
  \frac{d}{dc}h_{\theta}(c)=\frac{H_{\theta}^{\frac{2}{d}-1}(1-\alpha^2)c^{2- \alpha}}{B^3(c)}\left[\frac{(1-c^2)B(c)H_{\theta}(c)}{c^2}+\frac{(1-\alpha)}{d}\mathbb{E}_{\theta}(Ze^{-\frac{(1-\alpha^2)}{4B(c)}Z})\right]\,.
  \end{equation}
Focussing on the sign of the above expression, we have for $c>1$ and $T(c)=\frac{(c^2-1)B(c)}{c^2}$
$$  \frac{(1-c^2)B(c)H_{\theta}(c)}{c^2}+\frac{(1-\alpha)}{d}\mathbb{E}_{\theta}(Ze^{-\frac{(1-\alpha^2)}{4B(c)}Z})  > \frac{(1-\alpha) \epsilon}{d} - T(c)\,,$$
since $B(c)>2$ and $H_{\theta}(c)<1$.
Finally, since $\epsilon > 0 $ by assumption and since $T(c)$ increases in $c$, expression (\ref{hderivative}) is indeed positive for $\theta \in C$ whenever  
$ 1 < c \leq T^{-1}\left(\frac{(1-\alpha)\,\epsilon}{d}\right)= (k(d, \alpha, \sigma^2_X, \sigma^2_Y))^{1/2}$.   \\

For part {\bf (b)}, given the boundedness assumptions on the first and second moments of $T=Z$, it follows from Lemma \ref{bill} that $\epsilon \geq \underline{\epsilon} \, = b_0 \, e^{- \frac{1-\alpha^2}{8} \frac{b_2}{b_1} }\, > 0\,. $ The given lower bound on the cut-off point $k(d,\alpha, \sigma^2_X, \sigma^2_Y)$ follows with the lower bound $\underline{\tau}$ for $\tau$ and since $k(d,\alpha, \sigma^2_X, \sigma^2_Y)$ is increases as $\tau$ increases.  \qed 

\begin{remark}
\label{remarkafterTheorem}
\begin{itemize}

\item  Theorem \ref{main} holds as stated for spherically symmetric model $X \sim f_0(\|x-\theta\|^2)$ and $Y \sim N_d(\theta, \sigma^2_Y I_d)$ with known $f_0$.  The given proof applies throughout with the observation that Lemma \ref{alphadual} is an attribute of the normal density assumption on $Y$ only.

\item  Theorem \ref{main} represents a unified finding with respect to: {\bf (i)} the loss $L_{\alpha}$ for $|\alpha|<1$, {\bf (ii)} the choice of the plugged-in estimator $\hat{\theta}(X)$, {\bf (iii)} the dimension $d$, and {\bf (iv)} the parameter space $C$.  Moreover, the proof is unified.  As mentioned in the Introduction, the result adds to Fourdrinier et al. (2011)'s finding for Kullback-Leibler divergence loss.  Interestingly,  taking $\alpha=-1$ in part {\bf (a)} of Theorem \ref{main} leads to the cut-off point 
$$k(d,-1, \sigma^2_X, \sigma^2_Y) = 1 + \frac{\inf_{\theta \in C} 
\mathbb{E} (\|\hat{\theta}(X)-\theta\|^2}{d \sigma^2_Y})\,,$$  
which matches a sufficient condition given by Fourdrinier et al. (2011) and confirms further unification.   We expand further in Subsection 2.4 on the behaviour of the cut-off points $k(d,\alpha, \sigma^2_X, \sigma^2_Y)$.

\item  On the other hand,  Theorem \ref{main} does not provide a necessary and sufficient condition, as is the case for Theorems \ref{ax} and \ref{x+}, as well as the Kullback-Leibler finding of Fourdrinier et al. (2011).  An illustration below in Subsection 2.5 will further address this issue.  

\end{itemize}
\end{remark}

\begin{remark}
\label{constraints}
A large number of estimators $\hat{\theta}$ satisfy either the condition $\epsilon>0$   of part {\bf (a)} of Theorem \ref{main}, or part {\bf (b)}'s  boundedness conditions 
for $\mathbb{E}_{\theta} \|\hat{\theta}(X)-\theta\|^2$ and $\mathbb{E}_{\theta} \|\hat{\theta}(X)-\theta\|^4$.  Since $\mathbb{E}_{\theta} (Z \, e^{- (1-\alpha^2) Z/8})>0$ for all $\theta \in C$, the compactness of $C$ will suffice for the condition $\epsilon>0$ to be satisfied.  
For restricted but unbounded parameter spaces, and specifically for polyhedral cones $C$, which include orthant restrictions on some or all of the $\theta_i$'s, order constraints of the form $\theta_1 \leq  \cdots \leq \theta_d$, tree-order restrictions with $\theta_i \leq \theta_1$ for $i=2, \ldots, d$, and umbrella order restrictions of the form $\theta_1 \leq \cdots \theta_m \geq \theta_{m+1)} \geq \cdots \geq \theta_q$, and others, it follows from Marchand and Strawderman (2012) that $X$ is minimax under loss $\|\hat{\theta}-\theta\|^4$ with finite and constant minimax risk given by  $\mathbb{E}_{\theta} \|X-\theta\|^4 \, = \mathbb{E}_{\theta} (\|X-\theta\|^2)^2 \,=  (d^2 + 2d) \sigma^4_X$.  Consequently, estimators $\hat{
\theta}(X)$ that dominate $X$, such as projections onto $C$, will satisfy the conditions of Theorem \ref{main}.  We refer to Marchand and Strawderman (2012) for details and a list of further references.

Otherwise, we point out the following:

\begin{enumerate}
\item[{\bf (i)}]   The existence of a value $b_0>0$ arising with the lower bound condition on  $\mathbb{E}_{\theta} \|\hat{\theta}(X)-\theta\|^2$ is guaranteed with the condition that
$\hat{\theta}(X)$ be non-degenerate. 

\item[{\bf (ii)}]   The existence of $b_1 \in (0,\infty)$, related to the upper bound condition on  $\mathbb{E}_{\theta} \|\hat{\theta}(X)-\theta\|^2$, will be satisfied, for instance, by estimators $\hat{\theta}(X)$ that are minimax.  Many such choices are available in $d\geq 3$ dimensions or more.  On the other hand,  estimators with unbounded squared error loss, such as affine linear estimators $\hat{\theta}_a(X)=
aX$ studied in Section 2.1 will not satisfy the conditions.   Moreover,  the corresponding value of $\epsilon$ can be shown to be equal to $0$ for all $\sigma^2_X, \sigma^2_Y$, $\hat{\theta}_a$ with $0<a<1$, making Theorem \ref{main} inapplicable for such cases.  Of course, the analysis provided by Theorem \ref{ax} is    stronger anyway. 

\item[{\bf (iii)}]  For the existence of $b_2 \in (0,\infty)$ such that $\mathbb{E}_{\theta} \|\hat{\theta}(X)-\theta\|^4 \leq b_2 \sigma^4_Y$, it will suffice that $\hat{\theta}(X)$ dominate $X$ as a point estimator of $\theta$ under loss 
$\|\hat{\theta}-\theta\|^4$, as $\mathbb{E}_{\theta} \|X-\theta\|^4 \, = \, (d^2 + 2d) \sigma^4_X$, and one can thus choose $b_2= (d^2+ 2d) r^2$.  For $d \geq 3$, such estimators were obtained by Berger (1978) and include (for $\sigma^2_X=1$) James-Stein type estimators of the form $\hat{\theta}(X)= (1- \frac{f_1}{f_2+\|X\|^2}) X$ with $0 < f_1 \leq d-2$ and  $f_2 \geq 2 + \frac{d+1}{d+2} f_1$.  \\

\noindent Moreover, one can directly verify that Baranchik-type estimators (Baranchik, 1970) of the form $\hat{\theta}_{s(\cdot)}(X) = \left(1 - \frac{s(\|X\|^2)}{\|X\|^2}  \right) X$, with $s(\cdot) \geq 0$ have bounded $\mathbb{E}_{\theta}(\|\hat{\theta}_{s(\cdot)}(X) - \theta \|^4), \theta \in \mathbb{R}^d\,$, $d \geq 3$, as long as both $s(t)$ and $s(t)/t$ are bounded for $t>0$.  Such estimators include the positive-part James-Stein estimator $\hat{\theta}_{JS+}$ obtained with the choice $s(t) \, = \, \sigma^2_X \max\{t,d-2\}$. 

\item[{\bf (iv)}]    In related work for the model $X \sim N_d(\theta, \sigma^4_X I_d)$,  Fourdrinier, Ouassou and Strawderman (2008) provide various point estimators $\hat{\theta}$ that dominate $X$ under quartic loss
$\sum_{i} (\hat{\theta}_i - \theta_i)^4$.  These include James-Stein type estimators for $d \geq 5$ and Baranchik-type estimators (Baranchik, 1970) $\hat{\theta}_{s(\cdot)}(X)$ for $d \geq 3$ and certain conditions on $s(\cdot)$.   Specifically, their estimators $\hat{\theta}$ are such that:  
$$ \mathbb{E}\left(\sum_i (\hat{\theta}_i(X) - \theta_i)^4  \right) \leq  \mathbb{E}\left(\sum_i (X_i - \theta_i)^4  \right) \, = \, 3 d \sigma^4_X\,.$$ 
Lemma \ref{appendix} which follows, permits us to use such boundedness, coupled boundedness of quadratic loss, to guarantee that Theorem \ref{main} can be applied.

\begin{lemma} 
\label{appendix}
If $\hat{\theta}(X) \in \mathbb{R}^d$ is an estimator of $\theta \in \mathbb{R}^d$ such that 
$$\mathbb{E}\left(\sum_i (\hat{\theta}_i(X) - \theta_i)^4  \right) \leq  M_1\,, \hbox{ and } \mathbb{E}\left(\|\hat{\theta}(X)-\theta\|^2  \right) \leq M_2\,, \hbox{ for all } \theta\,,$$
then we have $\mathbb{E}\left(\|\hat{\theta}(X)-\theta\|^4 \right) \leq M_2^2 + dM_1$ for all $\theta$.
\end{lemma}
{\bf Proof.}  See Appendix.

It thus follows that estimators $\hat{\theta}(X)$ that dominate $X$ under quartic loss, that are also minimax under quadratic loss, are such that 
$\mathbb{E}\left(\|\hat{\theta}(X)-\theta\|^4 \right) \leq 4 d^2 \sigma^4_X$ for all $\theta \in \mathbb{R}^d$ by making use of Lemma \ref{appendix} with $M_1= 3 d \sigma^4_X$ and $M_2=d \sigma^4_X$.  For such estimators $\hat{\theta}$, Theorem \ref{main}'s $\epsilon$ is greater than $0$ and thus applies.

\end{enumerate}
\end{remark}

We conclude this subsection by pointing out that the dominating predictive density improvements that arise as a consequence of the above theorems, which are normal densities with an expanded variance, can be mixed to generate many other scale mixture of normals predictive densities which dominate the targeted plug-in density.  This is a consequence of Jensen's inequality, as laid out by the following.\footnote{A more general result appears in Yadegari, I. (2017). Prédiction, inférence sélective et quelques problèmes connexes.  Ph.D. thesis.  Université de Sherbrooke (http://savoirs.usherbrooke.ca/handle/11143/10167). }

\begin{lemma}
\label{Mixtureestimators}
Under the assumptions of Theorem \ref{main}, suppose that the predictive density $q_{\hat{\theta},c} \sim N_d(\hat{\theta(X)}, c^2 \sigma^2_Y I_d)$ dominates $q_{\hat{\theta},1}$ for $c \in (1,k^*)$ and under loss $L_{\alpha}$.  Let $F$ be a cdf such that $F(1)=0$ and $F(k^*)=1$. Then, the mixture density
$$\hat{q}_{F}\left(y; X \right)= \int_{1}^{k^*} q_{\hat{\theta},z}\left(y; X\right) \,dF(z) $$
also dominates $q_{\hat{\theta},1}$ under loss $L_{\alpha}$.
\end{lemma}
{\bf Proof.}
By Jensen's inequality, since $h_{\alpha}$ in (\ref{halpha}) is convex, we have with a change in order of integration
\begin{eqnarray*}
R_{\alpha}\left(\theta, \hat{q}_{F}\right) &=&  \mathbb{E} \left(\int_{\mathbb{R}^{d}} 
h_{\alpha}\left( \int_{1}^{k^*} \frac{\hat{q}_{z}\left(y; X\right)}{q\left(y \left| \theta \right. \right)} \,dF(z) \right)
q\left(y \left| \theta \right. \right)
\,dy \right)
\\
& \leq &  
\mathbb{E} \left( \int_{1}^{k^*} \int_{\mathbb{R}^{d}} 
 h_{\alpha}\left( \frac{\hat{q}_{z}\left(y; X\right)}{q\left(y \left| \theta \right. \right)} \right)
q\left(y \left| \theta \right. \right) \,dF(z)
\,dy \right).
\end{eqnarray*}
Now, use the assumed dominance results to infer that 
\begin{eqnarray*}
R_{\alpha}\left(\theta, \hat{q}_{F}\right) 
& \leq &  
\mathbb{E} \left( \int_{1}^{k^*} \int_{\mathbb{R}^{d}} 
 h_{\alpha}\left( \frac{\hat{q}\left(y; X\right)}{q\left(y \left| \theta \right. \right)} \right)
q\left(y \left| \theta \right. \right) \,dF(z)
\,dy \right)= R_{\alpha}\left(\theta, \hat{q}\right)
\end{eqnarray*}
 with strict inequality for at least one $\theta$, thus establishing the result. \qed

\subsection{Behaviour of the cut-off point $k(d,\alpha,\sigma^2_X,\sigma^2_Y)$ and simultaneous dominance}

We further expand here on the behaviour of Theorem \ref{main}'s cut-off point and implications for simultaneous dominance with respect to several losses $L_{\alpha}$.  In accordance with earlier examples, it seems plausible that $k(d,\alpha,\sigma^2_X,\sigma^2_Y)$ decreases in $\alpha$, $\alpha \in [-1,1)$. As an illustration, Figure 4 represents the cut-off points $k(d,\alpha, \sigma^2_X, \sigma^2_Y)$ for $d=3$, $\sigma^2_X=1$, $\sigma^2_Y=1,2,4$, $\alpha \in (-1,1)$, and the James-Stein estimator $\hat{\theta}_{JS}(X)= (1- \frac{d-2}{X'X}) X$.

\begin{figure}[!ht]
  \centering
    \includegraphics[width=0.99\textwidth]{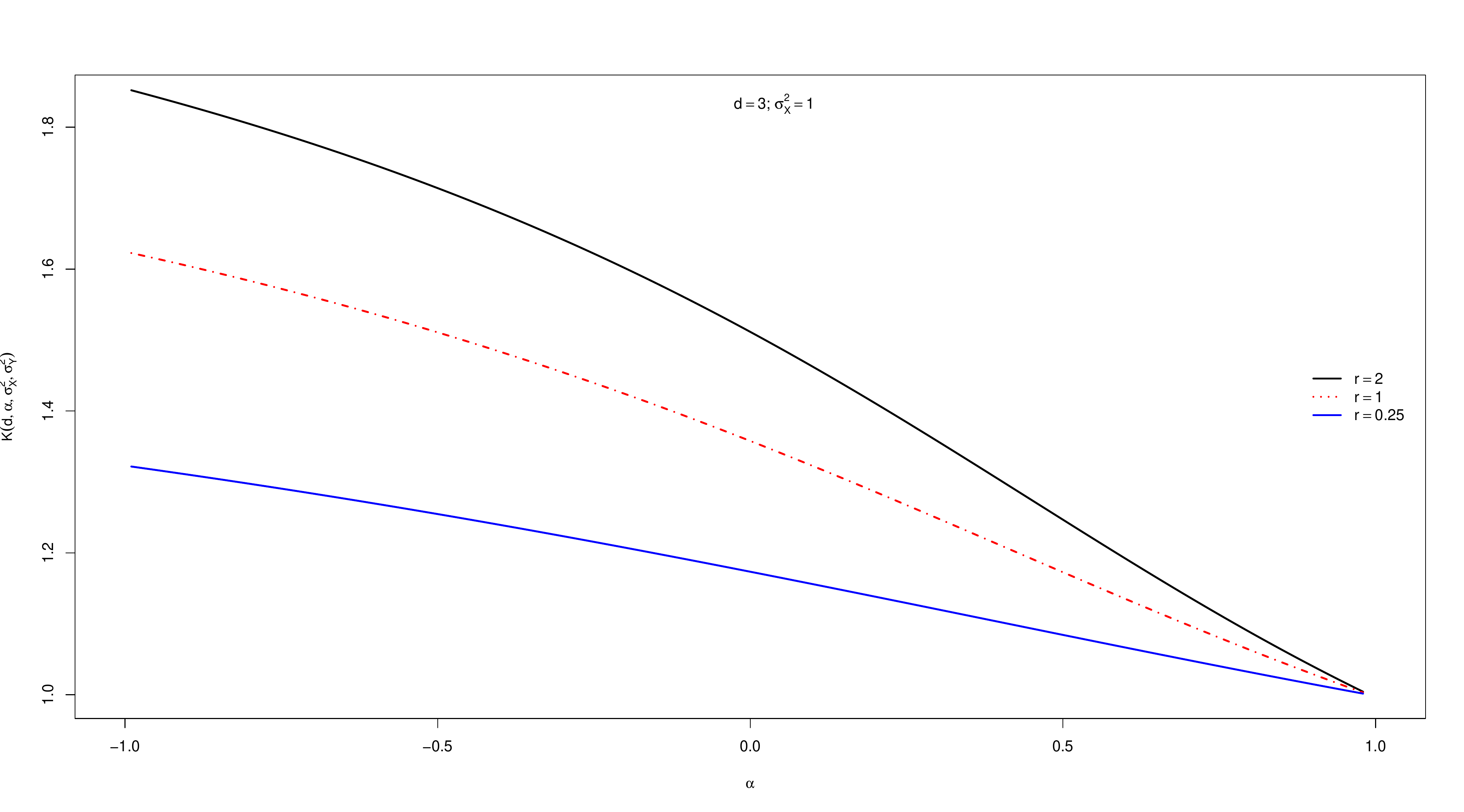}
     \label{kjs}
  \caption{Cut-off points $k(3,\alpha,1, \sigma^2_Y=1/r)$ as functions of $\alpha$, for the James-Stein estimator}
\end{figure}

It is of particular interest to focus on the benchmark Kullback-Leibler cut-off point for dominance.  Fourdrinier et al. (2011) show, for non-degenerate $\hat{\theta}$ other than $\hat{\theta}(X)=X$, that $q_{\hat{\theta},c}$ dominates the plug-in density 
$q_{\hat{\theta},1}$ for $1<c \leq 1+\underline{R}$ (as in Remark \ref{remarkafterTheorem}), and if and only if $1<c \leq c_0(1+\underline{R})$, with $$\underline{R}=  \frac{\inf_{\theta \in C} \mathbb{E} (\|\hat{\theta}(X)-\theta\|^2}{d \sigma^2_Y}\,,$$ and $c_0(t)$ the solution in 
$c \in (t,\infty)$ of the equation $(1-\frac{1}{c})t - \log c \, = \, 0$.    We pursue with an ordering between Hellinger and Kullback-Leibler cut-off points, as well as a monotonicity property, with implications for simultaneous dominance with respect to $L_{\alpha}$'s stated in the Corollary that follows.  The second part implies the first, but the alternative route for the first proof merits exposition.

\begin{theorem}
\label{alphadecreasing}
Consider the context of Theorem \ref{main} with a given $\hat{\theta}$ such that $\epsilon>0$, fixed $d, \sigma^2_X, \sigma^2_Y$ and consider the cut-off points $k(d,\alpha,\sigma^2_X, \sigma^2_Y)$.  Then, we have
\begin{enumerate}
\item[ {\bf (a)}]   $k(d,-1,\sigma^2_X,\sigma^2_Y) \geq k(d,0,\sigma^2_X,\sigma^2_Y)$;

\item[ {\bf (b)}] $k(d,\alpha, \sigma^2_X, \sigma^2_Y)$ is a non-increasing function of $\alpha \in [-1,0]$.
\end{enumerate}

\end{theorem}
{\bf Proof.}  See the Appendix for part {\bf (b)}.  For part {\bf (a)}, setting $\alpha=0$, we have $\tau=\epsilon(0) = d \underline{R}$, and 
$$ k(d,0,\sigma^2_X,\sigma^2_Y) \, = \, \frac{\underline{R}}{2} + \sqrt{(\frac{\underline{R}}{2})^2 + 1} \; \leq 1 + \underline{R}\, = k(d,-1,\sigma^2_X,\sigma^2_Y), $$
as $x+ \sqrt{x^2+1} \leq 1+2x$ for all $x>0$. \qed
 
As a consequence of the above, the following simultaneous dominance result is immediate.

\begin{corollary}
\label{simul}
Consider the context of Theorem \ref{main} with loss $L_{\alpha_0}$, $\alpha_0 \leq 0$, a given $\hat{\theta}$ such that $\epsilon>0$, and a value of $c^2 \in (1,k(d,\alpha_0,\sigma^2_X,\sigma^2_Y)$.   Then, the predictive density $q_{\hat{\theta},c}$ dominates the plug-in density $q_{\hat{\theta},1}$ for Kullback-Leibler loss as well as other losses $L_{\alpha}$ with $\alpha \leq \alpha_0$.
\end{corollary}
\subsection{Example}

\begin{example}

As seen above, Theorem \ref{main} is quite general and applies to many situations and many choices of the plug-in estimator $\hat{\theta}(X)$.  As an illustration, we focus on the positive-part James-Stein estimator given by $\hat{\theta}_{JS+}(X)$  (see part (iii) of Remark \ref{constraints}).

Theorem \ref{main} applies to the predictive densities  $q_{\hat{\theta}_{JS+},c}$ for $d \geq 3$, $\alpha \in (-1,1)$, $\sigma^2_X, \sigma^2_Y>0$, but we focus for the illustration on the roles of $d$ and $c$, and set $\alpha=0, \sigma^2_X=\sigma^2_Y=1$.  For implementing Theorem \ref{main}, we evaluate $\epsilon$ numerically, which yields $\tau$ and the cut-off point $k(d,0,1,1)$.  For $d=3$, we obtain for instance $\epsilon \approx 1.2009$ and thus $k(3,0,1,1) \approx 1.2200$.

Figure 5 compares the plug-in density with the variance expansion matching the cut-off points 
$k(d,0,1,1)$ for $d=3,5,7,9$.  The gains are moderate, or minimal, depending on $\|\theta\|$, and are 
further attenuated for larger dimension $d$ in accordance with the phenomenon exhibited at the outset of this section for the plug-in $\hat{\theta}(X)=X$.  

In opposition to earlier results, Theorem \ref{main} condition on the degree of variance expansion is not necessary and sufficient, so there is in theory room for improvement.   We proceeded with a numerical evaluation for $d=3$ giving dominance if and only if $1<c^2 \leq k^*$ with $k^* \approx 1.4883$.   An otherwise possible choice is given by the expansion  $c^2=1+(1-\alpha)r/2$, which is optimal for $\hat{\theta}(X)=X$, and equal to $c^2=3/2$ in our case.  For $d=3$, in accordance with the numerical evaluation, this does not lead to dominance, although gains are noticeable and more significant on a large part of the parameter space.     However, further numerical illustrations suggest worsened performance for larger $d$.   Finally, as a consequence of Corollary \ref{simul} we point out that the dominance illustrated here with the Hellinger cut-off points will hold for Kullback-Leibler loss, as well as all other $\alpha$-divergence choices with $\alpha \in (-1,0)$.

\begin{figure}[!ht]
  \centering
    \includegraphics[width=0.99\textwidth]{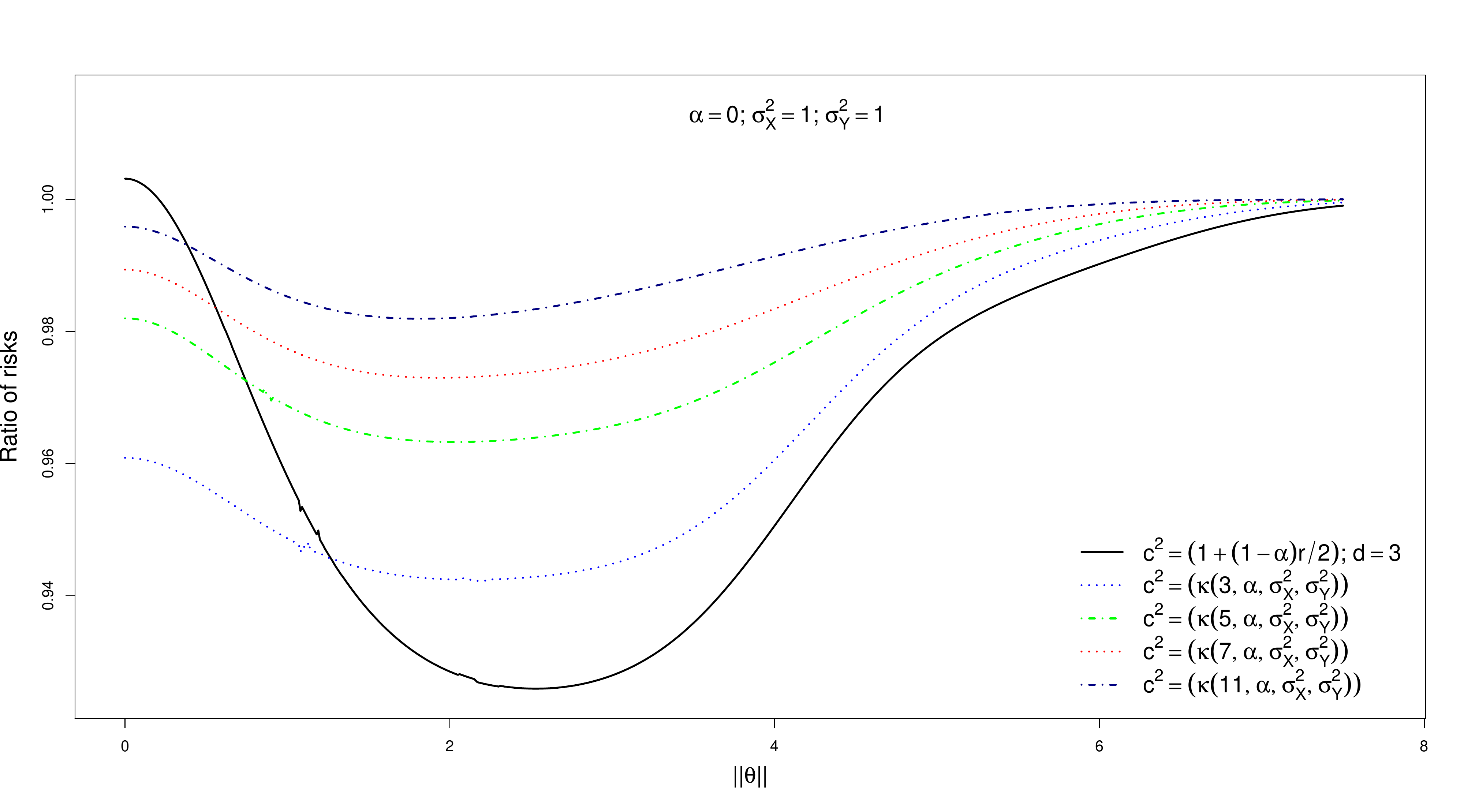}
     \label{ab2}
  \caption{Risk ratios $\frac{R_{\alpha}(\theta, q_{\hat{\theta}_{JS+},c})}{R_{\alpha}(\theta, q_{\hat{\theta}_{JS+},1})}$ for $c^2=k(d,\alpha=0,\sigma^2_X=1, \sigma^2_Y=1)$ as functions of $\|\theta\|$}
\end{figure}
\end{example}

\section{Concluding remarks}

For estimating the density of $Y \sim N_d(\theta, \sigma^2_Y I_d)$ based on $X \sim N_d(\theta, \sigma^2_X I_d)$, we establish the inadmissibility of plug-in densities $q_{\hat{\theta},1} \sim N_d(\hat{\theta}, \sigma^2_Y I_d)$ with respect to $\alpha-$divergence frequentist risk, and provide explicit variance expansion improvements of the form $q_{\hat{\theta},c} \sim N_d(\hat{\theta}, c^2 \sigma^2_Y I_d)$ with $c \in (1,c_0]$.  The results are quite general and apply to a large class of plugged-in estimators $\hat{\theta}(X)$.  Various implications arise, such as those with respect to robustness, as well as simultaneous dominance attained for a class of $\alpha-$divergence loss functions, including Kullback-Leibler.   Numerical illustrations complement the theory and are quite useful for instance in assessing the degree of improvement.

The findings are applicable in linear models, as well in the presence of normally distributed, or approximatively normally distributed, summary statistics that arise through sufficiency or in asymptotic settings.  The theoretical results in this paper highlight deficiencies present in the performance of plug-in densities and give credence to strategies to use alternatives.  It would be of interest, for instance, to develop Bayesian improvements and we feel the results here may serve such an objective.

\section{Appendix}

\subsection*{Proof of Lemma \ref{appendix}}
Let $\hat{\theta}=(\hat{\theta}_1, \ldots, \hat{\theta}_d)'$, $\theta=(\theta_1, \ldots, \theta_d)'$, and $Z_i=(\hat{\theta}_i - \theta_i)^2$ for $i=1, \ldots, d$.  Using the inequality $\mathbb{E}(Z_i Z_j) \leq \sqrt{\mathbb{E}(Z_i^2) \mathbb{E}(Z_j^2)} \, + \, \mathbb{E}(Z_i) \mathbb{E}(Z_j)$, we have
\begin{eqnarray*} \mathbb{E}\left(\|\hat{\theta}(X)-\theta\|^4 \right) & = & 
\sum_{i,j} \mathbb{E} (Z_i Z_j) \\
 & \leq &
\sum_{i,j} \sqrt{\mathbb{E}(Z_i^2) \mathbb{E}(Z_j^2)}  \, + \, \sum_{i,j} \mathbb{E}(Z_i) \mathbb{E}(Z_j) \\
& = &  d^2 \, \sum_{i,j} \frac{\sqrt{\mathbb{E}(Z_i^2) \mathbb{E}(Z_j^2)}}{d^2} + \left\lbrace\mathbb{E}\left(\|\hat{\theta}(X)-\theta\|^2 \right) \right\rbrace^2 \\
& \leq &  d^2 \sqrt{\frac{\sum_{i,j} \mathbb{E}(Z_i^2) \mathbb{E}(Z_j^2)}{d^2}} \, + \, M_2^2 \\
& = & d M_1 \; + \; M_2^2\,,
\end{eqnarray*}
the second inequality a consequence of Jensen's inequality applied to $g(t) = \sqrt{t}$ on $(0,\infty)$.  \qed

\subsection*{Proof of part (b) of Theorem \ref{alphadecreasing}}

From the definitions of $k(d,\alpha,\sigma^2_X,\sigma^2_Y)$, $\tau$, and $\epsilon$, 
we may write 
\begin{equation}
\label{k-alpha_1,alpha_2}
k(d,\alpha, \sigma^2_X, \sigma^2_Y) \, = \, \left( \beta(\alpha_1, \alpha_2) \, + \sqrt{\beta^2(\alpha_1, \alpha_2) + \frac{1+\alpha_1}{1-\alpha_1}}\, \right)_{\alpha_1=\alpha_2=\alpha}\,,
\end{equation}
with $\beta(\alpha_1, \alpha_2) \, = \frac{\epsilon(\alpha_2)}{2d} \, - \, \frac{\alpha_1}{1-\alpha_1}$ for $\alpha_1, \alpha_2 \in [-1,0]$.   With: {\bf (i)}  $\beta(\alpha_1, \alpha_2) \geq 0$ since $\epsilon(\alpha_2) >0$ and $\alpha_1 \leq 0$, {\bf (ii)} $\epsilon(\alpha_2)$ is non-increasing in $\alpha_2 \in [-1,0]$, {\bf (iii)} $\beta + \sqrt{\beta^2 + (1+\alpha_1)/(1-\alpha_1})$ is non-decreasing in $\beta \geq 0$, it follows that expression (\ref{k-alpha_1,alpha_2}) is, for fixed $\alpha_1 \in [-1,0]$ non-increasing in $\alpha_2 \in [-1,0]$.

To continue, it will thus suffice to show that (\ref{k-alpha_1,alpha_2}) is, for fixed $\alpha_2 \in [-1,0]$, non-increasing in $\alpha_1 \in [-1,0]$.  To this end, write (\ref{k-alpha_1,alpha_2}) as
\begin{equation}
\nonumber
T_{\alpha_2}(\alpha_1) \, = \, \frac{\epsilon(\alpha_2)}{2d} - \frac{\alpha_1}{1-\alpha_1} \, + \, \sqrt{\frac{\epsilon^2(\alpha_2)}{4d^2} + w(\alpha_1)}\,,
\end{equation}
with $ w(\alpha_1) \, = \, (\frac{\alpha_1}{1-\alpha_1})^2 \, - \, \frac{\alpha_1 \, \epsilon(\alpha_2)}{d(1-\alpha_1)} \, + \, \frac{1+\alpha_1}{1-\alpha_1}\,; \alpha_1 \in (-1,0]$.
Finally, we obtain
\begin{eqnarray*}
\frac{\partial}{\partial \alpha_1} T_{\alpha_2}(\alpha_1) \, & = & \, - \frac{1}{(1-\alpha_1)^2} \, + \, \frac{w'(\alpha_1)}{2 \sqrt{\frac{\epsilon^2(\alpha_2)}{4d^2} + w(\alpha_1)}} \\  & \leq &
\, \frac{1}{(1-\alpha_1)^3} \left(-(1-\alpha_1) \, + \, \frac{1}{\sqrt{\frac{\epsilon^2(\alpha_2)}{4d^2} + w(\alpha_1)}}  \right) \\
& \leq & 0\,,
\end{eqnarray*}
since $w'(\alpha_1) \leq \frac{\partial}{\partial \alpha_1} \left((\frac{\alpha_1}{1-\alpha_1})^2 \, + \, \frac{1+\alpha_1}{1-\alpha_1}  \right) = \frac{2}{(1-\alpha_1)^3}\,$, and 
$w(\alpha_1) \geq (\frac{\alpha_1}{1-\alpha_1})^2 \, + \, \frac{1+\alpha_1}{1-\alpha_1} \, = \frac{1}{(1-\alpha_1)^2} \geq 1$, for $\alpha_1 \in [-1,0]$. \qed

\section*{Acknowledgements}

Eric Marchand's research is supported in part by the Natural Sciences and Engineering Research Council of Canada.   We thank Bill Strawderman who provided the lower bound in (\ref{billA}).   Finally, we are grateful to Othmane Kortbi, Iraj Yadegari and Nasser Sadeghkhani for useful discussions on predictive density estimation under $\alpha-$divergence.


\end{document}